\documentclass[11pt]{article}

\usepackage{fullwidth}
\usepackage{color}
\usepackage[color,matrix,arrow]{xy}
\usepackage[top=1.5in, bottom=1.5in, left=1in, right=1in]{geometry}
\usepackage{amsgen}
\usepackage{amsmath}
\usepackage{amstext}
\usepackage{amsbsy}
\usepackage{amsopn}
\usepackage{amsfonts}
\usepackage{amssymb}
\usepackage{eepic}
\usepackage{graphicx}
\usepackage{tikz}
\usepackage{epsf}
\usepackage{pstricks}
\xyoption{all}

\def\Box{\square}

\def\tra#1{\smash{\mathop{\mid\kern
-1pt\joinrel\relbar\joinrel\relbar}\limits^{*}_{#1}}}
\def\longtra#1{\smash{\mathop{\mid\kern
-1pt\joinrel\relbar\joinrel\relbar\joinrel\relbar}\limits^{*}_{#1}}}
\def\vlongtra#1{\smash{\mathop{\mid\kern
-1pt\joinrel\relbar\joinrel\relbar\joinrel\relbar\joinrel\relbar}\limits^{*}_{#1}}}
\def\vvlongtra#1{\smash{\mathop{\mid\kern
-1pt\joinrel\relbar\joinrel\relbar\joinrel\relbar\joinrel\relbar\joinrel\relbar}\limits^{*}_{#1}}}
\def\vvvlongtra#1{\smash{\mathop{\mid\kern
-1pt\joinrel\relbar\joinrel\relbar\joinrel\relbar\joinrel\relbar\joinrel\relbar\joinrel\relbar}\limits^{*}_{#1}}}
\def\etra#1{\smash{\mathop{\mid\kern
-1pt\joinrel\relbar\joinrel\relbar}\limits_{#1}}}

\def\L{{\cal{L}}}
\def\M{{\cal{M}}}

\def\p{\varphi}

\def\bi{\begin{itemize}}
\def\ei{\end{itemize}}
\def\beq{\begin{equation}}
\def\eeq{\end{equation}}

\def\lf{\left\lfloor}   
\def\rf{\right\rfloor}

\newtheorem{T}{Theorem}[section]
\newcommand{\bt}{\begin{T}}
\newcommand{\et}{\end{T}}
\newcommand{\ftd}{$\square$\end{T}}

\newtheorem{Proposition}[T]{Proposition}
\newcommand{\bp}{\begin{Proposition}}
\newcommand{\ep}{\end{Proposition}}
\newcommand{\fpd}{$\square$\end{Proposition}}

\newtheorem{Lemma}[T]{Lemma}
\newcommand{\bl}{\begin{Lemma}}
\newcommand{\el}{\end{Lemma}}
\newcommand{\fld}{$\square$\end{Lemma}}

\newtheorem{Corol}[T]{Corollary}
\newcommand{\bc}{\begin{Corol}}
\newcommand{\ec}{\end{Corol}}
\newcommand{\fcd}{$\square$\end{Corol}}

\newtheorem{Result}[T]{Result}
\newcommand{\br}{\begin{Result}}
\newcommand{\er}{\end{Result}}
\newcommand{\frd}{$\square$\end{Result}}

\newtheorem{Example}[T]{Example}
\newcommand{\be}{\begin{Example}}
\newcommand{\ee}{\end{Example}}

\newtheorem{Problem}[T]{Problem}
\newcommand{\bq}{\begin{Problem}}
\newcommand{\eq}{\end{Problem}}

\newtheorem{Remark}[T]{Remark}
\newcommand{\brm}{\begin{Remark}}
\newcommand{\erm}{\end{Remark}}

\newtheorem{Definition}[T]{Definition}
\newcommand{\bd}{\begin{Definition}}
\newcommand{\ed}{\end{Definition}}

\newcommand{\proof}
   {\par\medbreak\noindent{\bf Proof}.\enspace}

\newcommand{\qed}{
$\Box$
\par\bigbreak}

\normalbaselines
\def\abstract#1{\par\bigskip
\begingroup\small
\baselineskip=12truept
\begin{center}ABSTRACT\end{center}
\par\medskip\par\noindent
\null\hfill\hbox{\vbox{\hsize=5truein\noindent#1}}
\hfill\null\par\endgroup\par}

\begin{document}
\title{Degree 2 Transformation Semigroups as Continuous Maps on Graphs: Foundations and Structure}
\author{{\bf Stuart Margolis, John Rhodes}}

\date{\today}

\maketitle

\begin{center}\small

\bigskip

\end{center}

\abstract{We develop the theory of transformation semigroups that have degree 2, that is, act by  partial functions on a finite set such that the inverse image of points have at most two elements. We show that the graph of fibers of such an action gives a deep connection between semigroup theory and graph theory. It is known that the Krohn-Rhodes complexity of a degree 2 action is at most 2. We show that the monoid of continuous maps on a graph is the translational hull of an appropriate 0-simple semigroup. We show how group mapping semigroups can be considered as regular covers of their right letter mapping image and relate this to their graph of fibers.} 


\section{Introduction}

A transformation semigroup (ts) is a pair $X=(Q,S)$, where $Q$ is a finite set and $S$ is a subsemigroup of the monoid $PF_{R}(Q)$ of all partial functions 
on $Q$, where composition is from left to right. A transformation monoid (tm) is a transformation semigroup $X=(Q,S)$ such that the identity function on $Q$ is in $S$.  Equivalently, a ts is a faithful right action of a semigroup $S$ on a finite set $Q$ by partial functions and a tm is a faithful right action of a monoid on $Q$. 
If $f:Q \rightarrow Q$ is a partial function, then $ker(f)$ is the equivalence relation on the domain of $f$ that identifies two elements $x,y$ if $xf = yf$. We call an equivalence class of $ker(f)$ a fiber of $f$. The degree, $deg(f)$ of $f$ is the maximal cardinality of a fiber and the degree, $deg(X)$ of a ts $X=(Q,S)$ is the maximum of the degrees of all elements of $S$. We write $Xd$ for the degree of a ts.

Thus a ts has degree 0 if and only if its semigroup consists of the empty function. We assume from hereon that every ts has a non-empty transformation, so that its degree is a positive integer. A ts $X=(Q,S)$ has degree 1 if and only if $S$ is a subsemigroup of the symmetric inverse monoid $I(Q)$ on $Q$, consisting of all partial bijections on $Q$. If deg(X)=1 and $S$ contains the inverse partial function $f^{-1}$ for every $f \in S$, then $S$ is an inverse semigroup \cite{Lawson, PetrichInv}. The theory of representation of inverse semigroups by partial 1-1 functions, that is, by degree 1 ts, is a central part of the theory of inverse semigroups and its applications \cite[Chapter 4]{PetrichInv}. 

It is straightforward to see that if $f,g \in PF_{R}(Q)$, then $deg(fg) \leq deg(f)deg(g)$. Also the collection of all partial functions of degree $k$ form a subsemigroup of $PF_{R}(Q)$ if and only if $k \leq 1$. Therefore studying ts of degree greater than 1 involves restrictions on the functions that appear in its semigroup. A result in \cite{Tilsonnumber} is that these restrictions can be encoded by the incidence structure whose blocks are the subsets of the fibers of the ts. The ts acts on the left by inverse image on this set and this embeds the semigroup $S$ into the translational hull \cite{qtheory} of the 0-simple semigroup whose structure matrix is the incidence matrix of the aforementioned incidence structure. The connection between the algebra of the ts and the geometry and combinatorics of the incidence structure gives a rich interaction between semigroup theory and incidence structures \cite{bibdtrans, bibdcont, projcont, AmigoWilson}.

The main purpose of this paper is to initiate a study of ts of a fixed degree $k > 1$ and we concentrate in this paper on ts of degree 2. In this case, the fibers, being of size at most 2 have the structure of a graph and the semigroup acts as partial functions on the vertices, such that the inverse image of each edge is either empty, a point or an edge. This allows results from graph theory to be used to understand the ts. We describe the structure of degree 2 ts in detail. 


We recall the basics of the complexity theory of finite semigroups and ts here. See \cite{qtheory} for more details. The complexity of a finite semigroup $S$ 
is the least number of non-trivial groups
needed in order to represent $S$ as a homomorphic image of a subsemigroup of a
wreath product of groups and semigroups whose maximal subgroups are trivial. Such a decomposition is guaranteed by the Krohn-Rhodes Theorem. The complexity problem is to compute this minimal number. Of course, any specific decomposition of $S$ of this type gives an upper bound of the complexity of $S$. Since the search
space for all shortest decompositions is infinite, there is no a priori reason that the
complexity of $S$ is computable. We write $Xc$ for the complexity of a transformation semigroup $X=(Q.S)$.

As far as we know, the only papers that study finite ts from the point of view of their degree are \cite{Tilsonnumber, TilsonMargolis}. Both of these relate degree to the Krohn-Rhodes complexity of a ts \cite{qtheory}. The main result of \cite{Tilsonnumber} shows that $Xc \leq Xd$, that is the degree of a ts $X$ is an upper bound for the complexity of $X$. In particular, since a ts $X=(Q,S)$  has complexity 0 if and only if every subgroup of $S$ is trivial, the complexity problem for transformation semigroups of degree 2 comes down to distinguishing between degree 2 ts of complexity 1 and those of complexity 2. We will explicate the structure of semigroups of degree 2. We then apply the topology of graphs and  to study the topological relationship between group mapping semigroups and right letter mapping semigroups \cite[Section 4.6]{qtheory}. 

We remark that some examples of ts with degree 2 have sophisticated structure from the point of view of complexity theory. This begins with the  so called Tall Fork Semigroup \cite[Section 4.14]{qtheory}, which has been used as a counterexample to a number of conjectures in complexity theory. In the continuation 
of this paper \cite{deg2part2}, we go well beyond the Tall Fork by giving examples that illustrate deeper aspects of complexity theory of degree 2 transformation semigroups and discuss the decidability question for transformation semigroups of degree 2 using the tools developed in \cite{Trans}.
 
\section{Transformation Semigroups of Degree 2 as Continuous Maps on Graphs}

All structures in this paper will be finite. As mentioned in the introduction, one must have some restrictions on a collection of partial functions of degree at most 2 to generate a semigroup of degree 2. Indeed it is well known that the monoid of all partial functions $PF_{R}(X)$ is generated by its group of units, that is the symmetric group on $X$ and by its elements of
rank $|X|-1$ which are clearly of degree at most 2. In fact, $PF_{R}(X)$ is generated by the symmetric group plus one idempotent of degree 2 and another of degree 1, both of rank $|X|-1$. Thus adding one element of degree 2 to the symmetric inverse semigroup, which consists of all the partial functions of degree at most 1, gives enough power to generate all partial functions. Thus, on the one hand ts of degree 2 must have a restricted structure and on the other hand, we will provide sophisticated examples of such ts.

We will show how to represent all transformation semigroups of degree 2 by the monoid of continuous partial functions in the sense of \cite{bibdtrans, bibdcont, AmigoWilson} on a graph. In this paper a graph is a pair $(V,E)$ where $V$ is a finite set and $E$ is a subset of the set of elements of cardinality of size 2. So these are simple graphs, there are no multiple edges nor loops. For ease of reading we write $uv$ if $\{u,v\} \in E$. When we need to look at structures with multiple edges  between vertices, we will explicitly use the term multigraph. 

We will also consider a graph $(V,E)$ as a simplicial complex of dimension one, whose faces are
$V \cup E$. In particular, we assume that each element of $V$ is a face. Recall that a simplicial complex is a pair $(V,H)$ where $V$ is a non-empty set and $H$ is a collection of non-empty subsets of $V$ such that if $Y \in H$ and $X\neq \emptyset \subseteq Y$, then $X \in H$.

Let $\Gamma=(V,E)$  and $\Gamma'=(V',E')$ be graphs. A partial function $f:V \rightarrow V'$ is {\em continuous} if for all $X \in V' \cup E'$, either $Xf^{-1}$ is empty or $Xf^{-1} \in V \cup E$. It is easy to prove that the composition of continuous functions is continuous as is the identity function. Thus we have a category
$\mathcal{C}\mathcal{G}$ whose objects are all finite graphs and whose morphisms are all continuous maps. 

We also define a stricter form of continuity of partial functions on graphs. Let $\Gamma =(V,E)$ and $\Gamma'=(V',E')$ be graphs with no isolated vertices.; A partial function $f:V \rightarrow V'$ is a {\em strict} continuous function if for all $uv \in E'$, either $(uv)f^{-1}$ is empty or $(uv)f^{-1} \in E$. That is, for strict partial functions, we do not allow $(uv)f^{-1}$ to be a an element of $V$. Strict functions, treat continuity in terms of the graph structure of $\Gamma$ whereas arbitrary continuity uses the simplicial complex structure of $\Gamma$. In an obvious way, we get the category $\mathcal{S}\mathcal{C}\mathcal{G}$ of strict continuous
functions on graphs. We discuss the connection between continuous and strict continuous functions below. 

We first contrast these notions with the usual definition of a morphism between graphs. 
Recall that a morphism from $\Gamma$ to $\Gamma'$ is a (total) function $f:V \rightarrow V'$
such that if $uv \in E$, then $(uf)(vf) \in E'$. The category $\mathcal{G}$ whose objects are all finite graphs and morphisms are all graph morphisms is one of the most extensively studied categories in mathematics. See for example, \cite{graphmorphs}.

Let $\Gamma=(\{1,2,3\},\{12,23\})$ be a path of length 2 and let $\Gamma'=(\{1,2\},\{12\})$ be a path of length 1.
Then the function $f:\{1,2,3\}\rightarrow \{1,2\}$ defined by $1f=1,2f=2,3f=1$ is a graph morphism, but it is not continuous since $(12)f^{-1}=\{1,2,3\}$ is not an edge.
Now consider the function  $g:\{1,2,3\}\rightarrow \{1,2,3\}$ defined by $1g=1,2g=1,3g=3$ is continuous on $\Gamma$, since $(12)g^{-1}=12$ and $(23)g^{-1}=3$, but is not a graph morphism since $(23)g=13$. Also, $g$ is continuous, but not strictly continuous.

The importance of continuous functions on graphs in our sense is that it gives us a way to describe all ts of degree at most 2. Let $M(\Gamma)$ denote the collection of all continuous partial functions from a graph $\Gamma$ to itself. This is the endomorphism monoid of $\Gamma$ in the category $\mathcal{C}\mathcal{G}$. The following  result shows that these give examples of ts of degree 2.

\bl \label{MGamma}

Let $\Gamma=(V,E)$ be a graph. Then $(V,M(\Gamma))$ is a tm of degree at most 2.

\el

\proof Clearly $M(\Gamma)$ is a submonoid of $PF_{R}(V)$. Therefore, $(V,M(\Gamma))$ is a tm. For all $v \in V$ and all $f \in M(\Gamma)$, if $vf^{-1}$ is non-empty, then $vf^{-1} \in V \cup E$ and thus has cardinality at most 2. \qed



We will give a number of examples of tm of the form $(V,M(\Gamma))$ for various graphs $\Gamma$, but first we prove a partial converse to Lemma \ref{MGamma}. We prove that every ts $X=(V,S)$ of degree at most 2 embeds into $(V,M(\Gamma(X)))$ where $\Gamma(X) = (V,E)$ is the graph whose edges are are the fibers of $X$ of size 2. We call $\Gamma(X)$ the graph of fibers of $X$. This is a special case of \cite[Theorem 3]{Tilsonnumber} which gives an analogous result for all degrees.

\bl \label{Embed}

Let $X = (V,S)$ be a ts of degree at most 2.  Then $S$ is a subsemigroup of 
$M(\Gamma(X))$ and thus $X$ embeds into $(V, M(\Gamma(X)))$.

\el 

\proof Let $X=(V,S)$ and let $\Gamma(X) =(V,E)$. If $f \in S$ and  $v \in V$, then $|vf^{-1}| \leq 2$ by assumption. In the case that $|vf^{-1}|=2$, then $vf^{-1} \in E$ by definition. If $xy \in E$, then there is a $w \in V$ and $g \in S$ such that $xy = wg^{-1}$. Then $(xy)f^{-1} = (wg^{-1})f^{-1}=w(fg)^{-1}$ is empty or    a fiber of $X$ and thus if non-empty belongs to $V \cup E$ by definition. Therefore $f \in M(\Gamma)$ and the Lemma is proved. \qed

Even though Lemma \ref{MGamma} and Lemma \ref{Embed} are straightforward to prove, they provide a bridge between semigroup theory and graph theory that we will exposit in the rest of the paper. 

We look a bit further at strict continuous functions. Let $\Gamma =(V,E)$ be a graph with no isolated vertices, so that for each $v \in V$, there is a $w \in V$ such that $vw \in E$. As above a partial function $f:V \rightarrow V$ is a {\em strict} continuous function if for all $uv \in E$, either $(uv)f^{-1}$ is empty or $(uv)f^{-1} \in E$. We let $SM(\Gamma)$ denote the collection of all strict continuous functions on $\Gamma$. Clearly, $SM(\Gamma)$ is a submonoid 
of $PF_{R}(V)$ and $(V,SM(\Gamma))$ is a tm. 

\bl \label{strict}

Let $\Gamma =(V,E)$ be a graph with no isolated vertices. Then $(V,SM(\Gamma))$ is a tm of degree 2. Furthermore, every strict continuous function is continuous.

\el

\proof Clearly $(V,SM(\Gamma))$ is a tm. We now show that its degree is 2.

Since $\Gamma$ has no isolated vertices, $|V| \geq 2$ and there is at least one edge $uv \in E$. Consider the partial function $f: V\rightarrow V$ with domain 
$\{u,v\}$ and such that $uf=vf=u$. Let $xy \in E$. If $u \in \{x,y\}$ then $(xy)f^{-1} = uv \in E$. Otherwise, $(xy)f^{-1} = \emptyset$. Therefore, $f \in S(\Gamma)$ and the degree of $(V,SM(\Gamma))$ is at least 2. 

Now let $f \in S(\Gamma)$. If $f$ is the empty function, then the degree of $f$ is 0. Otherwise, let $v \in Im(f)$. There is an edge $vw$ for some $w \in V$. Then 
$(vw)f^{-1} \in E$ and thus $|(vw)f^{-1}| = 2$. Therefore, $|vf^{-1}| \leq 2$ and thus $deg(f) \leq 2$ so the degree of $(V,SM(\Gamma))$ is 2.

Finally, if $f \in S(\Gamma)$, we know that if the inverse image of an edge of $\Gamma$ under $f$ is non-empty, then it is an edge and thus has cardinality 2. Since every vertex is contained in some edge, it follows that the inverse image of any vertex is size at most 2. Assume that the inverse image $|vf^{-1}| = 2$ so that $vf^{-1} = \{x,y\}$ for some 
$x \neq y \in V$. Let $w$ be such that $vw \in E$. Since the degree of $f$ is 2, it follows that $(vw)f^{-1} = \{x,y\}$ and thus $\{x,y\} \in E$. Therefore $f$ is continuous. \qed

As an example of a partial function that is continuous and  not strict continuous consider the complete graph $K_n$, where $n>1$. Then the partial function that is the identity on $\{1\}$ is continuous, but not strict continuous. We also note that the assumption that the graph $\Gamma$ has no isolated vertices is crucial in Lemma \ref{strict}. Consider the graph with no edges on a set $V$ with at least 2 elements. Then a partial function is continuous if and only if it is a partial bijection, whereas any partial function is strict continuous.

\subsection{The structure of continuous maps on graphs}

Let $\Gamma=(V,E)$ be a graph and $M = M(\Gamma)$ be the monoid of continuous functions on $\Gamma$. If $W \subseteq V$, let $\Gamma(W)$ be the subgraph induced by $W$, which is the graph with vertices $W$ and  all edges both of whose vertices belong to $W$. We first look at the structure of $\Gamma(Im(f))$, where $f \in M$.

The following lemma, showing that fibers are either vertices or edges is straightforward, but important.

\bl \label{sing1}

Let $f \in M$ and let $v \in Im(f)$ be such that $|vf^{-1}| = 2$. Then $vf^{-1} \in E$.

\el

\proof  By definition of continuous function,  $vf^{-1} \in V \cup E$ and thus must be in $E$ if  $|vf^{-1}| = 2$. \qed

Notice that $M$ is an ordered monoid with respect to restriction of partial functions. That is, if $f \in M$ and $g = f|_{Dom(g)}$, then $g \in M$, since if 
$X \in V \cup E$, then $Xg^{-1} \subseteq Xf^{-1}$. But $Xf^{-1} \in V \cup E$ and thus if $Xg^{-1}$ is either empty or belongs to $V \cup E$. This means by results of 
Stein \cite{Steinpartial} that the algebra of $M$ over a field is the same as a certain subcategory of the category of epimorphisms on $V$. This category is a so called $EI$-category. This allows for efficient calculation of the representation theory of $M$, a topic we will look at in a future paper.

We define $Sing(f)$, the singular part of $f$, to be the restriction of $f$ to the set
$\mathcal{M}(f) = \{v \in Dom(f)| \exists w \neq v, wf =vf\}$. That is, $\mathcal{M}(f)$ is the union of all fibers of $f$ of size 2. Recall that a matching in a graph is a collection of edges such that two distinct edges have no common vertex. An anti-clique, also called an independent set of vertices is a set of vertices
$W$ of $V$ such that $\Gamma(W)$ is the empty graph, that is, there are no edges between members of $W$.

\bl \label{sing2}

\begin{enumerate}

\item [(i)]{Let $f \in M$. Then the fibers of $Sing(f)$ form a matching of $\Gamma$ and the image of $Sing(f)$ is an anti-clique. 
Furthermore, $Sing(f)$ is strict continuous.}

\item [(ii)]{Conversely, let $\mathcal{M} =\{e_{1}, \ldots e_{m}\}$ be a matching of $\Gamma$ of size $m$ and let $X = \{x_{1}, \ldots , x_{m}\}$ be an 
anti-clique of size $m$ in $\Gamma$. Then the partial function $f$ that sends $e_{i}$ to $x_{i}, i=1, \ldots m$ is a strict continuous function and $f = Sing(f)$.}

\end{enumerate}

\el

\proof Let $f \in M$. Since $Sing(f)$ is a restriction of $f$, it follows that $Sing(f) \in M$. By definition, every fiber of $Sing(f)$ is of size 2 and thus is an edge of $\Gamma$ by Lemma \ref{sing1}. Thus, the fibers, forming a partition of the domain $\mathcal{M}(f)$ of $Sing(f)$ is a matching of $\Gamma$. If there were an edge $vw$ between two distinct vertices $v,w$ of the image of $Sing(f)$, then $|(vw)Sing(f)^{-1}|=4$, a contradiction. 

Let $uv \in E$. Since the image, $Im(Sing(f))$ of $Sing(f)$ is an anti-clique, 
$|uv \cap Im(Sing(f))| \leq 1$. Therefore, either $(uv)(Sing(f))^{-1}$ is the empty set or 
$(uv)(Sing(f))^{-1}$ is a fiber of $f$ and thus an edge of $\Gamma$ by definition. Therefore, $Sing(f)$ is a strict continuous function.  This proves item (i).

Let $M$, $X$ and $f$ be as in item (ii) and let $uv \in  E$. Then $|uv \cap X| \leq 1$ since $X$ is an anti-clique. Therefore, $(uv)f^{-1}$ is either the empty set or 
$(uv)f^{-1}$ is a fiber of $f$ and thus an edge of $\Gamma$ by definition. Therefore, $f$ is strict continuous. Clearly, $f = Sing(f)$. \qed

Let $f \in M$. We define the injective part of $f$ to be the restriction of $f$ to  $Dom(f) - \mathcal{M}(f)$ and denote this by $Inj(f)$. If $f$ and $g$ are partial functions on the same set with disjoint domains, then their join $f \vee g$ is the partial function whose domain is the union of the domains of $f$ and $g$ and such that $x(f \vee g)= xf$ if $x \in Dom(f)$ and $x(f \vee g) = xg$ if $x \in Dom(g)$. In particular, $f=Sing(f) \vee Inj(f)$ is the join decomposition of $f$ into its singular and injective parts. We now characterize which functions can serve as $Sing(f)$ and $Inj(f)$ for an arbitrary continuous function $f$. Let $k(f)=\{v \in Im(f)| \exists{w \in Im(f), vw \in E}\}$.  Thus $k(f)$ is the largest subgraph of $\Gamma(Im(f))$ with no isolated vertices.

\bl \label{sing}

\begin{enumerate}

	\item[(i)]{Let $f$ be a continuous partial function on a graph $\Gamma$. Then the partition $ker(Sing(f))$ is a matching of $\Gamma$ and the image $Im(Sing(f))$ is an anti-clique in $\Gamma$. Moreover, $ Im(Sing(f)) \cap Im(Inj(f)) = \emptyset$, the inverse partial function of $Inj(f)$ is a bijective graph morphism onto its image, and there are no edges between a vertex in $Im(Sing(f))$ and a vertex in $Im(Inj(f))$.}

 \item[(ii)]{Conversely, let $g$  be a partial function such that $ker(g)$ is a matching in $\Gamma$ and $Im(g)$ is an anti-clique in $\Gamma$. Let $h$ be a partial bijection such that $h^{-1}$ is a bijective graph morphism onto its image. Assume further that $Dom(g) \cap Dom(h) = Im(g) \cap Im(h) = \emptyset$ and that there are no edges between $Im(g)$ and $Im(h)$. Then $f=g \vee h$ is continuous and $g=Sing(f), h=Inj(f)$.} 

\end{enumerate} 
\el

\proof By Lemma \ref{sing2}, $ker(Sing(f))$ is a matching and $Im(Sing(f))$ is an anti-clique in $\Gamma$. Assume that $v \in Im(Sing(f)) \cap Im(Inj(f))$. Since 
$Dom(Sing(f)) \cap Dom(Inj(f)) = \emptyset$, it follows that $|vf^{-1}|=3$ contradicting the continuity of $f$. Therefore $Im(Sing(f)) \cap Im(Inj(f))=\emptyset$.
The subgraph $\Gamma(Inj(f))$ consists of the subgraph $k(f)$ of $\Gamma$ together with the collection of isolated vertices $Im(Inj(f)) \setminus  k(f)$.  We now claim that $Inj(f)^{-1}$ is a bijective graph morphism onto its image. Indeed, $Inj(f)^{-1}$ is a bijection onto its image and if $uv$ is an edge of $Im(Inj(f))$, then $(uv)Inj(f)^{-1}$ is an edge by the definition of continuity. This proves item (i).

Now let $g$ and $h$ be as in item (ii) and let $f = g \vee h$. Let $v \in V$. If $vf^{-1}$ is not empty then either $v \in Im(h)$ and $|vf^{-1}| = 1$ or $v \in Im(g)$, 
whence $vf^{-1} \in E$, since each fiber is an edge in a matching of $\Gamma$. Assume that $uv \in E$ and that $(uv)f^{-1}$ is non-empty. If $uv \subseteq Im(h)$, 
then $(uv)f^{-1}=(uv)h^{-1} \in E$, since $h^{-1}$ is a bijective graph morphism onto its image. Otherwise, the assumptions imply that $uv$ intersects exactly one of $Im(g)$ or 
$Im(h)$ in exactly one vertex. If $(uv) \cap Im(h)$ is a vertex, then $|(uv)|f^{-1}=1$ since $h$ is injective. If $(uv) \cap Im(g)$ is a vertex, then $(uv)f^{-1} \in E$, since each fiber of $g$ is an edge. Therefore, $f$ is continuous. Finally, since each fiber of $g$ has exactly two elements and each fiber of $h$ has exactly one element, it follows that $g=Sing(f)$ and $h=Inj(f)$. \qed









Before continuing, a few examples are useful. We will look at the structures of these monoids later in the paper. 

\be \textbf{Empty Graphs}

Let $V$ be a non-empty set and let $N(V)$ be the graph with no edges on $V$. Then a function is continuous if and only if it is a partial bijection, that is, $|vf^{-1}| = 1$
for all $v \in Im(f)$. Therefore $M(N(V))$ is the symmetric inverse monoid on $V$ consisting of all partial bijections on $V$.

\ee

\be \textbf{Complete Graphs} \label{completeexample}

Let $n>0$ and let $K_n$ be the complete graph on $n$ vertices $V_{n}=\{1, \ldots n\}$.  Let $f$ be a continuous function on $K_n$. By Lemma \ref{sing2}, the image of 
$Sing(f)$ is an anti-clique and thus $|Im(Sing(f))| \leq 1$. If $Im(Sing(f))$ is empty, then $f$ is a partial bijection on $V_{n}$. Clearly any partial bijection is continuous on $K_n$. Assume then that $Im(Sing(f)) = \{v\}$ for some $v \in V_{n}$. Then by Lemma \ref{sing} for all $w \in Im(Inj(f))$, $vw$ is not an edge.  Therefore, $Inj(f)$ is the empty function and $f$ is a partial constant function with domain an edge $e$. Let 
$e \in E, v \in V$. Define $f_{e,v}:V \rightarrow V$ be the partial constant function with domain $e$ that sends both vertices of $e$ to $v$. Then $f_{e,v}$ is continuous by the above. Conversely, any partial constant function with domain two elements is continuous. Therefore $M(K_{n})$ consists of the symmetric inverse monoid on $V$ together with the collection $\{f_{e,v}| e \in E, v \in V\}$.


\ee

The following Lemma will determine strict continuous functions on connected graphs.

\bl  \label{strinj}

Let $\Gamma = (V,E)$ be a connected graph with $|V| > 1$ and let $f$ be a strict continuous non-empty partial function on $\Gamma$. Suppose that there is a $v \in V$ such that
$|vf^{-1}|=1$. Then $f$ is an automorphism of $\Gamma$. 

\el

\proof Let $v \in V$ be such that $|vf^{-1}|=1$ and let $w$ be a neighbor of
$v$, which exists since $\Gamma$ is connected. $(vw)f^{-1} \in E$, since it is non-empty and $f$ is strict continuous. It follows that $|wf^{-1}|=1$. Therefore, $f^{-1}$ is defined on every vertex in the neighborhood of $v$ and the inverse image of each neighbor of $v$ has cardinality 1. By a straightforward induction on the distance of a vertex $x$ from $v$ in $\Gamma$, $f^{-1}$ is defined on all of $V$ and is a bijection. Clearly a continuous bijection of $\Gamma$ is an automorphism of $\Gamma$ and we are done. \qed 

The following corollary follows immediately from this Lemma.

\bc \label{strconn}

Let $\Gamma = (V,E)$ be a connected graph with $|V| > 1$. Then the monoid $SM(\Gamma)$ of all strict continuous functions on $\Gamma$ consists of the empty function, the group of automorphisms and all continuous functions on $\Gamma$ of the form $f =Sing(f)$. The latter are precisely the continuous functions such that every fiber has exactly two elements.

\ec

We note by Example \ref{completeexample} that every partial bijection of $V_n$ is continuous on the complete graph $K_n$. It follows from Lemma \ref{strinj} that 
if $n > 1$ that there are injective continuous functions on $K_n$ that are not strict. It is easy to see that the strict continuous functions of $K_n$ consist of the empty function, the partial constant functions $f_{e,v}$ that are defined on the edge $e$ and sends both elements of $e$ to $v \in V$ together with the symmetric group 
$S_n$ which is the automorphism group of $K_n$.

Let $\Gamma$ be a non-connected graph with no isolated vertices and let $C$ be a connected component of $\Gamma$. Clearly, the identity function restricted to $C$ is a non-empty strict continuous function that is not an automorphism of $\Gamma$. Thus the assumption that $\Gamma$ is connected is necessary in Lemma \ref{strinj}. An easy modification of Lemma \ref{strinj} shows that the domain of a strict continuous injective function on an arbitrary graph is a union of connected components.

\be

Let $\Gamma = K_{n,n}$, the complete bipartite graph on $2n$ vertices, with bipartition $V_{n} = B_{n} \dot{\cup} W_{n}$ where $B_{n}=\{b_{1} \ldots, b_{n}\}$ and
$W_{n} = \{w_{1}, \ldots w_{n}\}$. We compute the monoid of strict continuous functions $SM(K_{n,n})$. By Corollary \ref{strconn}, we need to compute the singular functions, that is, those continuous functions of the form $f=Sing(f)$ and all automorphisms of $K_{n,n}$. It is well known that the automorphism group of $K_{n,n}$ is the group of the permutation group $Z_{2} \wr (\{1,\ldots, n\},Sym(n))$ where $Z_2$ is the group of order 2 and $Sym(n)$ is the symmetric group on $n$ elements. Thus the group of units has $2^{n}n!$ elements.

A singular function on $K_{n,n}$ has a matching as its domain. It is clear that a matching in $K_{n,n}$ can be identified with a partial bijection 
$M:B_{n} \rightarrow W_{n}$. The range of a singular function is an anti-clique and thus lies wholly in either $B_{n}$ or $W_{n}$. Thus a singular continuous function of rank $k$ is given by a partial bijection $M:B_{n} \rightarrow W_{n}$ with a set of of cardinality $k$ as domain. The corresponding matching is the set of all edges
of the form $bM(b)$, where $b \in Dom(M)$. There are $(\binom{n}{k})^{2}k!$ partial bijections of rank $k$, since we can choose any $k$-set of $B_{n}$ as domain and any $k$-set of $W_{n}$ as range and we can permute the range as we please for each such choice of set in $W_{n}$. For each such partial bijection $M$ we have a unique singular function on $K_{n,n}$ by using the associated matching of $M$ as fibers and sending them arbitrarily to either a $k$-set in $B_{n}$ or a $k$-set in $W_{n}$. Thus each such $M$ give $2(\binom{n}{k})k!$ singular functions, since we can permute the image at will. Therefore there are $2(\binom{n}{k})^{3}(k!)^{2}$ singular functions of rank $k$. Call this number $S_k$. Therefore the number of singular functions is $\sum_{k=0}^{n}S_{k}$.

\ee

Thus, $SM(K_{n,n})$ is a rather large semigroup from the point of view of cardinality. It also has a rich ideal structure as will be discussed in the next section. Despite this, the complexity of every monoid of continuous functions on any graph has complexity at most 2, as we will now show.

\section{On the complexity of degree 2 transformation semigroups} \label{compdeg2}

The main result of \cite{Tilsonnumber} is that the degree of a transformation semigroup is an upper bound to its complexity. We adapt the proof in that paper to the case of transformation semigroups of degree 2. 

We recall the basics of the complexity theory of finite semigroups and ts here. See \cite{qtheory} for more details. The complexity of a finite semigroup $S$ 
is the least number of non-trivial groups
needed in order to represent $S$ as a homomorphic image of a subsemigroup of a
wreath product of groups and semigroups whose maximal subgroups are trivial. Such a decomposition is guaranteed by the Krohn-Rhodes Theorem. The complexity problem is to compute this minimal number. Of course, any specific decomposition of $S$ of this type gives an upper bound of the complexity of $S$. Since the search
space for all shortest decompositions is infinite, there is no a priori reason that the
complexity of $S$ is computable. We write $Xc$ for the complexity of a transformation semigroup $X=(Q.S)$.

We need some further concepts from the complexity theory of finite transformation semigroups. See \cite{qtheory} for more details. In this paper, a relation
$\p:U\rightarrow V$ between sets $U$ and $V$ is a function $U\rightarrow 2^{V}$, from $U$ to the power set $2^{V}$ of $V$. The graph of $\p$ denoted $\#\p$ is the usual subset of $U \times V, \#\p=\{(u,v)|v \in u\p\}$.

 Let $X=(Q,S)$ and $Y=(P,T)$ be transformation semigroups. Following \cite[Section 4.14]{qtheory}, a relation $\phi:Q \rightarrow P$ is called a 
{\em relational morphism} of transformation semigroups, written $X \triangleleft_{\phi} Y$ if 

\begin{enumerate}

\item [(1)]{$\phi$ is fully defined. That is, for all $q \in Q, q\p \neq \emptyset$.}
	
\item [(2)]{ For every $s \in S$, there is an element $\overline{s} \in T$ such that for all $p \in P$, $(p\phi^{-1})s \subseteq (p\overline{s})\phi^{-1}$. That is,
for all $q \in Q, s \in S$ such that $qs$ is defined and for all $p \in P$ with $p \in q\phi$, $y\overline{s}$ is defined and belongs to $xs\p$. In this case we say that $\overline{s}$
covers $s$.}

\end{enumerate}

 We remark that relational morphisms of ts are the inverse of what are called relational covers in \cite{Eilenberg}. We use the notion of relational morphism as we are following the presentation in Chapter 4 of \cite{qtheory}.

A relational morphism is called a {\em division} if $\phi^{-1}$ is a partial function. In this case we write $X \prec Y$. Clearly if $X \prec Y$, then $Xc \leq Yc$. 
The companion relation of $\p$ is the relation $f_{\p}:S \rightarrow T$ defined by $sf_{\p}=\{t \in T| t \text{ covers } s\}$. Recall that a relational morphism between semigroups is a fully defined relation $f:S \rightarrow T$ such that $(s_{1})f(s_{2})f \subseteq (s_{1}s_{2})f$ for all $s_{1},s_{2} \in S$. Equivalently, $\#f$ is
a subsemigroup of $S \times T$ that projects onto $S$. It is easy to prove that if $\p$ is a relational morphism of ts, then $f_{\p}:S \rightarrow T$ is a relational
morphism of semigroups. A {\em parametrized relational morphism} between ts $\Phi:(Q,S)\rightarrow Y=(P,T)$ is a pair $\Phi=(\p_{1},\p_{2})$ consisting of a relational morphism of ts
$\p_{1}:X \rightarrow Y$ and a relational morphism of semigroups $\p_{2}:S \rightarrow T$, such that $\#\p_{2}$ is contained in the graph of the 
companion relation $\#f_{\p_{1}}$. The canonical paratmetrization of a relational morphism $\p_{1}:X \rightarrow Y$ is the pair $\Phi=(\phi_{1},f_{\p_{1}})$.


By an automaton we mean a pair $\mathcal{A}=(Q,\Sigma)$ where $Q$ (the states) and $\Sigma$ (the alphabet) are finite sets together with a partial function 
$\delta:Q \times \Sigma \rightarrow Q$. As usual, we denote the image of $(q,x)$ by $qx$ for all $q \in Q, x \in \Sigma$ such that $\delta(q,x)$ is defined. This is 
the usual notion of a deterministic partial automaton without initial nor terminal states. Thus each $x \in \Sigma$ defines a partial function from $Q$ to itself. The 
transformation semigroup defined by $\mathcal{A}$ is the ts $(Q,S)$, where $S$ is the subsemigroup of the monoid of partial functions on $Q$ generated by 
all $x \in \Sigma$.

The derived automaton $\mathcal{A}_{\Phi}=(Q_{\Phi},\Sigma_{\Phi})$ of a relational morphism $\Phi:(Q,S) \rightarrow (P,T)$ parametrized by $(\phi_{1},\p_{2})$ is the automaton
$D_{\Phi} = (Q_{\Phi},\Sigma_{\Phi})$ where $Q_{\Phi}=\{(q,p)|p \in q\phi_{1}\}$, the graph of $\phi_{1}$ and 
$\Sigma_{\Phi} = \{(p,(s,t)) \in P \times \#\p_{2} | pt \neq \emptyset\}$. The action on $Q_{\Phi}$ is given by defining:


\begin{equation*}
    (q,p)(p',(s,t)) = \begin{cases}
              (qs,pt)               &  p=p' \text{ and } qs \text{ is defined}\\
               \text{undefined}               & \text{otherwise}\\
           \end{cases}
\end{equation*}

The derived ts $\mathcal{D}_{\Phi} =(Q_{\Phi},S_{\Phi})$ is the transformation semigroup of $\mathcal{A}_{\Phi}$. The importance of the derived transformation semigroup is given by the following theorem of Tilson. See \cite[Chapter 3, Section 8]{Eilenberg} for a proof. We emphasize that opposed to this paper, which follows 
\cite[Section 4.14]{qtheory}, \cite{Eilenberg} uses relational covers of ts, which are the inverse relation of our relational morphism of ts. It is straightforward to translate from one of these notions to the other.

\bt \label{derived}

Let $X\triangleleft_{\Phi} Y$ be a parametrized relational morphism of transformation semigroup. Then $X \prec D_{\Phi} \wr Y$.

\et

Let $\Gamma = (V,E)$ be a graph and let $M(\Gamma)=(V,M(\Gamma))$ be its transformation monoid of continuous functions. Let $P$ be the set of anti-cliques of $\Gamma$.

\bl \label{action}

Let $Y$ be an anti-clique in $\Gamma$ and let $f \in M(\Gamma)$. Then $Y(Inj(f)) \cup Im(Sing(f))$ is an anti-clique.


\el

\proof Assume that $uv$ is an edge that belongs to  $(Y(Inj(f)) \cup Im(Sing(f)))$. If $uv \in Y(Inj(f))$, then $(uv)f^{-1} \subseteq Y$ is an edge by continuity of $f$. But $Y$ is an 
anti-clique, so this is not possible. We have seen in Lemma \ref{sing} that $Im(Sing(f))$ is an anti-clique, so $uv$ is not contained in $Im(Sing(f))$. Therefore, the only possibility left is that (without loss of generality) $u \in Im(Sing(f))$ and $v \in Im(Inj(f))$. But then $|(uv)f^{-1}|=3$ contradicting continuity. \qed

Now consider the symmetric inverse monoid $S=SIS(V)$. Then $S$ acts on the right of the power set $2^V$ by total functions by  sending $X \subseteq V$ to its image $Xs$
for $s \in S$. We can thus construct the semidirect product, $S \rtimes 2^V$ of $S$ with the semilattice of subsets $2^V$ under union by defining $(s,X)(t,Y)=(st,Xt \cup Y)$, where $s,t \in S$ and  $X,Y \subseteq V$. 

The next lemma requires the Fundamental Lemma of Complexity \cite[Section 4.9]{qtheory}. This states that if $R:S \rightarrow T$ is a surjective morphism, such that the inverse image of each idempotent of $T$ has only trivial subgroups, then the complexity of $S$ is equal to the complexity of $T$. Furthermore, it is known that the symmetric inverse monoid on $V$ has complexity 0 if $|V| \leq 1$ and has complexity 1 otherwise. See Chapter 4 of \cite{qtheory} for more details.

\bl 

The complexity of $S \rtimes 2^V$ is 0 if $|V| \leq 1$ and is 1 otherwise.

\el

\proof Consider the projection morphism $P: (S \rtimes 2^V) \rightarrow S$ sending $(s,X)$ to $s$. Idempotents in $S$ are restictions of the identity function to subsets of $V$. Let $ X \subseteq V$. Then $(1|_{X})P^{-1}=\{(1_X,Y)|Y \subseteq V\}$. A direct calculation shows that every element in this set is itself an idempotent.  Therefore, the Fundamental Lemma of Complexity implies that the complexity of $S \rtimes 2^V$ is the same as that of $S$ and the lemma is proved. \qed

\br

One can check that the inverse image of an idempotent of the projection morphism $P: (S \rtimes 2^V) \rightarrow S$ is actually an $\mathcal{L}$-trivial band. That is, every $\mathcal{L}$-class is trivial, or equivalently it satisfies the identity $yx = xyx$ along with $x^{2}=x$. This gives a tighter decomposition result, but is not necessary for the purposes of this paper.

\er

The semigroup $(S \rtimes 2^V)$ acts on the right of $2^V$ by defining $X(s,Y) = Xs \cup Y$ for $s \in S, X,Y \subseteq V$. Let $T'$ denote the subsemigroup of $(S \rtimes 2^V)$ generated by $\{(Inj(f),Im(Sing(f)))|f \in M(\Gamma)\}$. Then by Lemma \ref{action},  the subset $P=\{X \subseteq V|X \text{ is an anti-clique}\}$ is invariant under this action. Let $T$ be the faithful image of $T'$ under this action. Then we have the transformation semigroup $Y = (P,T)$.

\bt \label{contcomp}

Let $(Q,S)$ be a transformation semigroup of degree 2. Then the complexity of $(Q,S)$ is at most 2.

\et

\proof By Lemma \ref{Embed} it is enough to prove the assertion for $(V,\M(\Gamma))$, the transformation semigroup of continuous maps for some graph 
$\Gamma =(V,E)$. Consider the transformation semigroup $Y=(P,T)$ we defined above, where $P=\{X \subseteq V|X \text{ is an anti-clique}\}$. 

We define a relation $\phi_{1}: V \rightarrow P$ by sending $v \in V$ to $\{A| v \in A\}$. Since every singleton set is an anti-clique, $\phi_{1}$ is a fully 
defined relation. Let $s \in M(\Gamma)$. We claim that $s$ is covered by $(Inj(s),Im(Sing(s))) \in T$ (or more precisely, the element of $T$ represented by this pair). Indeed, if $A$ is an anti-clique, then $(A\phi_{1}^{-1})=\{v|v \in A\}=A$ so that $(A\phi_{1}^{-1})s=As=A(Inj(s)) \vee Sing(s)) \subseteq (A(Inj(s)) \cup Im(Sing(s))=(A(Inj(s)),Im(Sing(s)))\phi_{1}^{-1}$. Therefore, $(V,M(\Gamma)) \triangleleft_{\phi_{1}} (P,T)$ is a relational morphism of transformation semigroups. 

We let $\Phi$ be the canonical parametrization $\Phi=(\phi_{1},f_{\phi_{1}})$. Consider the derived transformation semigroup $D_{\Phi}$ relative to this parametrization. A non-empty function in $D_{\Phi}$ is of the form $(A,s,A(Inj(s)) \cup Im(Sing(s)))$, where $A$ is an anti-clique. We claim that $s$ acts as a partial one-to-one map on $A$. Indeed, assume that $u\neq v \in A$ and that 
$us=vs =x \in V$. Then by continuity it follows that $xs^{-1} =\{u,v\}$ is an edge, contradicting that $A$ is an anti-clique. Therefore, $D_{\Phi}$ divides the symmetric inverse monoid acting on the states of $D_{\Phi}$. As mentioned above, this implies that $D_{\Phi}c \leq 1$. Finally, by Theorem \ref{derived}, 
$(V,M(\Gamma)) \prec D_{\Phi} \wr (P,T)$. We have seen that both $D_{\Phi}c \leq 1$ and $(P,T)c \leq 1$ and therefore $(V,M(\Gamma))c \leq 2$ as desired. \qed

This theorem gives another big difference between the category $\mathcal{G}$ of graphs and graph morphisms and the category 
$\mathcal{C}\mathcal{G}$ of graphs and continuous functions.
The well known theorem of Frucht that every finite group is the automorphism group of a finite graph can be extended to prove that every finite monoid is the endomorphism morphism of a finite graph in $\mathcal{G}$ \cite[Theorem 1.35]{graphmorphs}. In contrast to this, it is known that there are monoids of arbitrary complexity and thus the collection of finite monoids that have faithful representations by ts of degree at most 2, that is the collection of submonoids of
endomorphism monoids in $\mathcal{C}\mathcal{G}$ is a proper collection of monoids. For example, the monoid of all partial functions on a finite 
set $X$ has complexity $|X|-1$ \cite[Chapter 4]{qtheory}. Theorem \ref{contcomp} also shows that any monoid of complexity strictly greater than 2, has no faithful representation by a ts of 
degree 2.

\be \label{small2}

We construct the smallest transformation semigroup of degree 2 and complexity 2. Let $\Gamma = (\{1,2,3,4\},\{12,23,34,41\})$ be a 4-cycle. There are two perfect matchings, namely, $12|34$ and $14|23$ and two maximal anti-cliques, $13$ and $24$. All 8 possibilities of assigning a perfect matching to a maximal anti-clique defines a singular continuous map. Thus, for example, sending $12 \mapsto 3$, $34\mapsto 1$ defines a continuous function and all 7 other possibilities do as well. We add the permutation with cycle decomposition $z = (12)(34)$ which is an automorphism of $\Gamma$ and thus an invertible continuous map. A straightforward calculation shows that the semigroup generated by the 8 singular maps of rank 2 and $z$ is a monoid of order 10. One can check that this monoid is precisely the 10 element monoid of \cite[Example 4.10.12]{qtheory}, which is known to be the unique semigroup of order 10 of complexity 2 and that all semigroups of order at most 9, have complexity at most 1.

\ee

As a corollary to the proof of Theorem \ref{contcomp} we get the following result. By a singular transformation semigroup of degree 2, we mean a ts $(Q,S)$ such that each $s \in S$ acts as a singular partial function on $Q$. That is, for all $q \in Im(s)$, $|qs^{-1}|=2$.

\bc \label{singcomp}

Let $\Gamma=(V,E)$ be a graph and let $X=(V,S)$ be a singular transformation semigroup of continuous functions on $\Gamma$. Then the complexity of $X$ is at most one.

\ec

\proof In the proof of Theorem \ref{contcomp}, every singular function $s$ is covered by the element $(\theta, Im(s)) \in (S \rtimes 2^V)$, where $\theta$ is the empty function on $V$. Since for all $X,Y \subseteq V$, we have $(\theta,X)(\theta,Y) = (\theta,Y)$, the semigroup $T$ in the relational morphism of ts constructed in the Theorem is a right-zero semigroup. Since we proved that $X \prec (D_{\phi}) \wr T$, it follows that the complexity of $X$ is at most 1. \qed

\bc \label{singcomp2}

Let $V=(G,E)$ be a connected graph with trivial automorphism group. Then the monoid $SM(\Gamma)$ of strict continuous functions has complexity at most 1.

\ec

\proof It follows from Corollary \ref{strconn} that $SM(\Gamma)$ consists of the singular continuous functions and the identity transformation. By the previous corollary, the singular part has complexity 1. Since adding an identity to a semigroup is easily seen to preserve its complexity, we are done. \qed

\brm

It follows from Example \ref{small2} that the monoid of strict continuous functions on a 4-cycle has complexity 2, so the assumption that the automorphism group of the graph is trivial is necessary in Corollary \ref{singcomp2}.

It is known that asymptotically, almost all graphs have trivial automorphism groups (and even trivial endomorphism monoids in the category $\mathcal{G}$ of graph morphisms) \cite[Theorem 4.7]{graphmorphs} and thus asymptotically almost all graphs have strict monoids of complexity at most 1.

\erm

Despite having complexity at most 2, the ideal structure of monoids of continuous and strict continuous maps can be intricate. If $S$ is a semigroup, define $S\delta$ to be the longest chain of regular $\mathcal{J}$-classes of $S$. It is well known that $Sc \leq S\delta$. In fact, the Depth Decomposition Theorem 
\cite[Theorem 4.9.15]{qtheory} shows that the longest chain of regular $\mathcal{J}$-classes that contain non-trivial groups is an upper bound to complexity, but we don't need that result for this discussion. It is a well known fact that $S\delta$ is equal to the longest chain of idempotents of $S$ in the usual idempotent order. Recall that if $e,f$ are idempotents of $S$, then $e \leq f$ if and only if $e=ef=fe$ gives this order. Also, it is easy to prove, since idempotents and  $\mathcal{J}$-classes lift under morphisms of finite semigroups that if $S \prec T$, then $S\delta \leq T\delta$. Therefore, if $(V,S)$ is a transformation semigroup, 
then $S\delta \leq (PT_R{V})\delta \leq (|V|+1)$. The next Lemma shows that the monoid $M(\Gamma)$ of continuous functions on a graph $\Gamma =(V,E)$ always attains this bound.

\bl \label{contdepth}

Let $\Gamma=(V,E)$ be a graph and let $M(\Gamma)$ be its monoid of continuous functions. Then $M(\Gamma)\delta = |V|+1$.

\el

\proof Let $V =\{v_{1} \ldots , v_{n}\}$. It is clear that for every subset $X$ of $V$ the identity restricted to $X$, $1_{X}$, is a continuous function. Therefore, taking $X_{i}=\{v_{1},\ldots,v_{i}\}, 0\leq i\leq |V|$, we have a chain of idempotents $1_{X_{0}} < \ldots < 1=1_{X_{|V|}}$ and the result is proved. \qed

For the monoid of strict continuous functions $SC(\Gamma)$, the situation is different. It follows  from the discussion following Corollary \ref{strconn} that for all $n>1$ the depth of $SM(K_{n})$ is 3 where $K_n$ is the complete graph on $n$ vertices. More generally, we have the following result. Recall that if $\Gamma=(V,E)$ is a graph, then its matching number $\Gamma\nu$ is the size of a maximal matching of $\Gamma$. The independence number of $\Gamma$,  $\Gamma\alpha$ is the size of a maximal anti-clique in 
$\Gamma$.

\bl

Let $\Gamma=(V,E)$ be a connected graph and let $SM(\Gamma)$ be its monoid of strict continuous function. Then $SM(\Gamma)\delta \leq \text{ min}\{\Gamma\nu,\Gamma\alpha\} + 2$.

\el

\proof By Corollary \ref{strconn} a strict continuous function on $\Gamma$ is either the empty function, an automorphism or a singular function. We have seen that the partition of every singular function is a matching of $\Gamma$ and the corresponding image is an anti-clique. The result follows. \qed

We see now why $SM(K_{n})\delta=3$, since $K_{n}\nu = \lf \frac{n}{2} \rf$ and $K_{n}\alpha=1$. On the other hand, for the complete bipartite graph $K_{n,n}$, we have 
$K_{n,n}\nu=K_{n,n}\alpha=n$ and thus $SM(K_{n,n})\delta = n+2$. We don't know if there are results which connect the matching number and independence number of a graph. It is known that the matching number of a graph $\Gamma$ is equal to the independence number of its line graph $\L(\Gamma)$. Recall that the line graph has the edges of $\Gamma$ as vertices and has an edge for each pair of edges of $\Gamma$ that have a vertex in common.

\section{Semilocal Theory and Translational Hulls} \label{semiloc}

We now give a connection between the monoids we have been studying with the classical notions of translational hulls of 0-simple semigroups and the semilocal theory of finite semigroups. Thus this section enters coordinates into the coordinate-free approach we've taken up to now. We refer the reader to  \cite[Chapter 4.6, Chapter 5.5]{qtheory} for background details on these ideas. We recall some of these notions and fix some notation.

A finite semigroup $S$, possibly with a 0-element is 0-simple, if $S^{2} \neq \{0\}$ and every two-sided ideal of $S$ is either $S$ or $\{0\}$. By the Rees Theorem \cite[Appendix A.4]{qtheory}, these are isomorphic to Regular Rees matrix semigroups $\mathcal{M}^{0}(G,A,B,C)$ where $G$ is a finite group, $A$ and $B$ are sets and $C:B \times A \rightarrow G^{0}$. Regularity means that for every $ b \in B$ there is an $a \in A$ with $C(b,a) \neq 0$ and dually, for every $a \in A$, there is a $b \in B$
with $C(b,a) \neq 0$. The underlying set is $(A \times G \times B) \cup \{0\}$ and multiplication is such that 0 is the zero-element and
$(a,g,b)(a',g',b') = (a,gC(b,a')g',b')$ if $C(b,a') \neq 0$ and to 0 otherwise.

We can view the function $C$ as a matrix with rows labeled by $B$ and columns labeled by $A$. Regularity then means that each row and column have at least one non-zero element. It will be convenient for us to sometimes use ``inner product'' notation. We write $<b,a>_{C} = C(b,a)$. 

We can view the monoid $PF_{R}(X) (PF_{L}(X))$ as the monoid of $|X| \times |X|$  row monomial (column monomial) matrices over $\{0,1\}$. Row monomial (column monomial) means that each row (column) of the matrix has at most one non-zero element. We identify $f \in PF_{R}(X)$ with the matrix $X_f$ with $X_{f}(x,xf)=1$ for 
each $x \in Dom(f)$ and 0 otherwise. More generally, if $G$ is a finite group, then the wreath product $G \wr PF_{R}(X)$ can be identified with the monoid 
of row monomial matrices over $G^{0}$ and the reverse wreath product $PF_{L}(X) \wr_{r} G$ with the monoid of column monomial matrices over $G$. More explicitly, if 
$(\p,f) \in G \wr PF_{R}(X)$, where $f \in PF_{R}(X), \p:Dom(f) \rightarrow G$ is identified with the row monomial matrix $X_{(\p,f)}$ where $X_{(\p,f)}(x,xf)=x\p$ if 
$x \in Dom(f)$ and 0 otherwise. A dual definition gives the connection between column monomial matrices over $G$ and the reverse wreath product $PF_{L}(X) \wr_{r} G$. See \cite[Section 5.5]{qtheory} for more details.

Given this identification, we define the translational hull $\Omega(\mathcal{M}^{0}(G,A,B,C))$ of a zero simple semigroup $\mathcal{M}^{0}(G,A,B,C)$ to be the set of pairs $(X,Y)$ where $X$ is a $|B| \times |B|$ row monomial matrix over $G$ (that is, a member of $G \wr PF_{R}(B)$) and $Y$ is a $|A| \times |A|$ column monomial matrix
over $G$ (that is, a member of $PF_{L}(X) \wr_{r} G$) such that $XC=CY$. It is easy to check that $\Omega(\mathcal{M}^{0}(G,A,B,C))$ is a monoid containing 
an isomorphic copy of $\mathcal{M}^{0}(G,A,B,C)$ as its unique 0-minimal ideal. An important property of $\mathcal{M}^{0}(G,A,B,C)$ is that if $S$ is any semigroup that has $\mathcal{M}^{0}(G,A,B,C)$ as its unique 0-minimal ideal and $S$ acts faithfully on either the right or the left of $\mathcal{M}^{0}(G,A,B,C)$, then there is an embedding of $S$ into $\Omega(\mathcal{M}^{0}(G,A,B,C))$ that is an isomorphism on $\mathcal{M}^{0}(G,A,B,C)$. 
 
We translate the matrix equation $XC=CY$ into wreath product notation. Then $X=X_{(\p,f)}$ for a unique element $(\p,f) \in G \wr PF_{R}(X)$ and $Y =Y_{(f^{*},\p^{*})}$ for a unique element $(f^{*},\p^{*}) \in PF_{L}(X) \wr_{r} G$. Direct matrix multiplication gives 
\begin{equation}\label{linked}  b\p<bf,a>_{C}=<b,f^{*}a>_{C}\p^{*}a \end{equation} for all $b \in B, a \in A$ in inner product notation. These are called the {\em linked equations}. 

Let's specialize the above to the case that $G=\{1\}$, the trivial group. Let $(X,Y) \in \Omega(\mathcal{M}^{0}(\{1\},A,B,C))$. Then $X$ is the row monomial matrix of a unique partial function $f \in PF_{R}(B)$ and $Y$ is the column monomial matrix of a unique partial function $f^{*} \in PF_{L}(A)$. Then $XC=CY$ translates into 
\begin{equation}\label{trivlinked}<bf,a>_{C}=<b,f^{*}a>_{C}\end{equation} for all $b \in B, a \in A$ in inner product notation. Thus $f$ and $f^{*}$ are adjoint with respect to $<,>_{C}$. We write $<,>:B \times A \rightarrow \{0,1\}$ instead of $<,>_C$ if $C$ is understood and fixed.

If $G$ is the trivial group, we can view $C$ as the incidence matrix of an incidence system with points $B$ and blocks given by considering the row of $a \in A$ as
the subset $\overline{a}=\{b \in B|<b,a>=1\}$. Then membership in the translational hull, $<bf,a>=<b,f^{*}a>$ means that for all $b \in B, a \in A$, $bf \in \overline{a}$ 
if and only if $b \in \overline{f^{*}(a)}$. That is, $\overline{f^{*}(a)}$ is the inverse image $f^{-1}(\overline{a})$. Thus the blocks of the incidence system corresponding to $C$ is
 closed under inverse image with respect to $f$. This connection clearly motivates our monoids of continuous and strict continuous monoids on a graph. We make this
connection precise below. Viewing the translational hull of 0-simple semigroups over the trivial group as ``continuous'' partial maps on the corresponding incidence system has proved to be a fruitful connection between semigroup theory and combinatorics \cite{bibdtrans, bibdcont, projcont, AmigoWilson}.

Let $\Gamma = (V,E)$ be a graph. The graph incidence matrix of $V$ is the $|V| \times |E|$ matrix $S=S(\Gamma)$ whose entry in position $(v,e)$ is 1 if $v$ is a vertex of $e$ and 0 otherwise. Note that every column of $S$ has exactly two non-zero entries. The row sum of row $v$ is the degree of $v$ in $\Gamma$, where here, degree is used as in graph theory as the number of edges on which $v$ is a vertex. Thus there is a row of all zeros if and only if there is an isolated vertex in $\Gamma$. We will assume when talking about the graph incidence matrix that the graph has no isolated vertices. Since we are working with simple graphs (no loops nor multiple edges) all columns of $S$ are distinct. Two rows $v,w$ of $S$ are the same if and only if $vw$ is a connected component of $\Gamma$.

If we view $\Gamma$ as a simplicial complex, we have the simplicial incidence matrix, which is the $|V| \times |V \cup E|$ matrix $C = C(\Gamma)$ with entries 
$C(v,w)=1$ if and only if $v=w, v,w \in V$ and as above, $C(v,e)=1$ if and only if $v$ is a vertex of $e$.  As matrices, the relationship between $C$ and $S$ is that 
$C = [S|I_{V}]$, where $I_{V}$ is the $|V| \times |V|$ identity matrix. That is, we add $|V|$ columns to $S$ which contain the identity matrix in order to build $C$ from $S$. It is clear then that distinct columns and rows of $C$ are not equal to one another.

We use these matrices as the structure matrices of the following Rees matrix semigroups over the trivial group: $\mathcal{M}^{0}(1,E,V,S)$ and 
$\mathcal{M}^{0}(1,V \cup E,V,C)$. Assuming that $\Gamma$ has no isolated vertices means that these are regular Rees matrix semigroups and thus 0-simple semigroups. 
The discussion above leads immediately to the following result.

\bl \label{transhull}

Let $\Gamma=(V,E)$ be a graph. Then the translational hull of $\mathcal{M}^{0}(1,V \cup E,V,C)$ is isomorphic to the monoid $M(\Gamma)$ of all continuous partial functions
on $\Gamma$. If $\Gamma$ has no isolated vertices, then the translational hull of $\mathcal{M}^{0}(1,E,V,S)$ is isomorphic to the monoid $SM(\Gamma)$ of strict continuous partial functions on $\Gamma$.

\el

\section{Group Mapping Semigroups as Regular Covering Spaces Over Right Letter Mapping Semigroups}

In this section we look at coverings of graphs and their relation to the semi-local theory of semigroups. See \cite[Chapter 4.6]{qtheory} and \cite[Chapters 7,8]{Arbib} for background material. Here we allow graphs with multiple edges and loops as is usual in the topology of graphs. To distinguish between simple graphs, that we've used up to now, we will use the term multigraph when wanting to emphasize that multiple edges and loops are possible. We note that the definition of continuous functions makes sense for multigraphs as well. 

By Lemma \ref{transhull} we can view both the monoid of continuous functions and the monoid of strict continuous functions of a graph as the translational hull of 0-simple semigroups over the trivial group. Recall \cite[Section 4.6]{qtheory} that a Group-Mapping (GM) semigroup $S$ has a unique 0-minimal minimal regular ideal with a non-trivial maximal subgroup on which $S$ acts faithfully on both the left and right. Equivalently, $S$ is a subsemigroup of the translational hull of a 0-simple semigroup $I(S)$ with non-trivial maximal subgroup, containing $I(S)$ and such that no rows or columns of the structure matrix of $I(S)$ are proportional.

A Right Letter Mapping (RLM) semigroup is a semigroup that acts faithfully on the right of a unique 0-minimal ideal, that has a trivial maximal subgroup. Every GM semigroup $S$ acts on the $\mathcal{L}$-classes of $I(S)$ by right multiplication by sending for $s \in S$ the $\mathcal{L}$-class $L$ to $Ls$, which is also an
$\mathcal{L}$-class of $S$. The image of this representation is an $RLM$ semigroup denoted by $RLM(S)$. An important fact is that the complexity of $S$ is at most one more than the complexity of RLM($S$). Furthermore the problem of deciding the complexity of an arbitrary finite semigroup can be reduced to deciding whether the complexity of a GM semigroup $S$ is equal to that of RLM($S$) or not. See Section 4.6 of \cite{qtheory} for details.

In this section we look at the graph of fibers of GM semigroups of degree 2 as regular covering spaces over the free action by the maximal subgroup of its 0-minimal ideal. In the second part of this paper \cite{deg2part2} we use this to describe examples of GM semigroups of degree 2 that need the sophisticated tools from \cite{Trans} to compute their complexity. The main tool are graphs whose edges are labeled by elements of a group $G$. These go by a surprisingly large number of names: gain graphs (sign graphs when the group is cyclic of order 2), voltage graphs \cite{topgraph}, $G$-labeled graph \cite[Chapter 4.13]{qtheory} and more. They are related to torsors and to 
principal $G$-bundles over a group $G$ as well \cite{MM-Sheaves}. We will use the term gain graph in this paper. There is a very extensive literature on these structures. See Thomas Zaslavsky's 550 page (and growing) bibliography on gain graphs and related structures obtainable at http://people.math.binghamton.edu/zaslav/Bsg/index.html.

We emphasize the connection between gain graphs, $GM$ semigroups, their $RLM$ image and wreath products. A $GM$ semigroup $S$ with maximal subgroup $G$ in $I(S)$ acts by continuous functions on the derived graph (defined below) of a gain graph over its maximal subgroup $G$. The derived graph turns out to be the graph of fibers of $S$. This gives a topological interpretation of the embedding of $S$ into 
$G \wr RLM(S)$ as embedding $S$ into the wreath product of $G$ and the monoid of continuous functions on the graph of fibers of $RLM(S)$. 

Let $\Gamma = (V,E)$ be a multigraph. Choose an orientation $e^{+}$ for each $e \in E$. If the orientation $e^{+}$ starts at $v$ and ends at $w$ for $e = vw$, then we define the opposite edge to be $e^{-}$ that starts at $w$ and ends at $v$. A more precise version of this would be to define graphs in the sense of Serre \cite{Serre} but this description suits our needs in this paper. There is an evident notion of oriented path in a graph. When we write path, we mean oriented path. Each path $p=e_{1}^{\epsilon_{1}} \ldots e_{n}^{\epsilon_{n}}$ has an opposite path 
$p^{-1}=e_{n}^{-\epsilon_{n}} \ldots e_{1}^{-\epsilon_{1}}$. Here, $\epsilon_{i} \in \{+,-\}$ and $-\epsilon$ is the opposite sign to that of $\epsilon$. The collection of all paths under
concatenation of paths and opposite is then the free category with involution $P(\Gamma)$ on $\Gamma$.

Two oriented paths $p_{1}$ and $p_{2}$ are elementary homotopic if $p_{2}$ is obtained from $p_{1}$ by insertion or deletion of an oriented path of the form 
$e^{+}e^{-}$ or $e^{-}e^{+}$  into $p_{1}$. Two paths are homotopic if there is a finite sequence of elementary homotopies starting with $p$ and ending with $q$. It is easy to see that homotopy is a congruence on the category $P(\Gamma)$. The quotient is easily seen to be a groupoid, a category in which each morphism is an isomorphism. This is the free groupoid $P_{1}(\Gamma)$ on $\Gamma$, sometimes called the fundamental groupoid of $\Gamma$. If $v \in V$, then the collection $P_{1}(\Gamma , v)$ of all (congruence classes ) of paths beginning and ending at $v$ is a group called the fundamental group of $\Gamma$ at $v$. It is well known \cite{LyndonandSchupp} that each element of $P_{1}(\Gamma)$ is represented by a unique reduced path, that is, a path having no path of length two of the form $e^{+}e^{-}$ or $e^{-}e^{+}$ and that $P_{1}(\Gamma , v)$ is a free group for all $v \in V$. All of this is a combinatorial version of the fact that the fundamental groupoid of the geometric realization of a  graph is a free groupoid and that the fundamental group at a vertex $v$ is a free group.

Let $G$ be a group. A $G$-labeling of a graph $\Gamma$ with an orientation is a map $l:E^{+} \rightarrow G$ from the set of positively oriented edges to $G$. We define
$l(e^{-})=l(e^{+})^{-1}$. Clearly we can extend $l$ to a map from $P_{1}(\Gamma)$ to $G$ and thus consider $l$ to be a functor from $P_{1}(\Gamma)$ to $G$, considered as a one-object groupoid. We call a graph with a $G$-labeling, $\Gamma(l)=(V,E,l)$ a gain graph over $G$. 

Let  $\Gamma(l)=(V,E,l)$ be a gain graph over $G$. The derived graph 
$\Delta((\Gamma(l))$ is the ordinary graph with vertex set $G \times V$ and edge set $G \times E$. If the directed edge $e^{+}$ runs from $v$ to $w$ in the graph $\Gamma$ and if $l(e^{+}) = h$, then the directed edge $(g,e)^{+}$ runs from $(g,v)$ to $(gh,w)$. 

Gain graphs are essentially the same as regular covering projections of graphs. They are usually called voltage graphs in this context. We recall the basics. See \cite[Chapters 1 and 2]{topgraph} for details. Let $\Gamma' = (V',E')$ be a graph. Let $G$ be a group that acts freely (that is, without fixed points) on the left by graph automorphisms on the vertices and edges of $\Gamma'$. Then the vertex and edge orbits of $G$ are easily seen to form a graph $\Gamma'/G$. The natural projection from 
$\Gamma'$ to $\Gamma'/G$ is called a regular covering projection. They are the combinatorial version of regular coverings of graphs in the topological sense. More generally, a graph morphism between graphs $\Gamma' = (V',E')$ and $\Gamma = (V,E)$ is called a regular covering if it is equivalent to a regular projection. See
\cite[Section 1.3.8]{topgraph} for details.

Let $\Gamma(l)$ be a gain graph over the group $G$ and graph $\Gamma$. Then $G$ acts freely on the vertices (edges) of the derived graph $\Delta(\Gamma(l))$ by 
$g(h,v)=(gh,v)$ ($g(h,e)=(gh,e)$), where $g,h \in G, v \in V, e \in E$. It is clear from the definitions that the quotient of $\Delta(\Gamma(l))$ by this action of $G$
is isomorphic to $\Gamma$ and is a regular covering projection. The following Theorem \cite[Theorem 2.2.2]{topgraph} of Gross and Tucker gives the converse of this result. Section \cite[Section 2.4]{topgraph} treats the case of arbitrary covering projections, but we do not need these results in this paper.
 
\bt \label{covergain}

Let $G$ be a group acting freely on the graph $\Gamma'$ and let $\Gamma$ be the resulting quotient graph. Then there is a gain graph $\Gamma(l)$ with labels from $G$ and underlying graph $\Gamma$ such that the derived graph of $\Gamma(l)$ is isomorphic to $\Gamma'$ and the projection from this derived graph to $\Gamma$ is equivalent to the regular covering projection from $\Gamma'$ to $\Gamma$.

\et


We now look at regular coverings in the case of graphs of fibers of transformation semigroups of degree 2. We then relate this to the projection from a GM semigroup to its RLM image. Let $X=(V,S)$ be a transformation semigroup of degree 2 and let $\Gamma = (V,E)$ be its graph of fibers. Let $G$ be a group and let 
$l:E^{+} \rightarrow G$ for some fixed orientation of the edges of $\Gamma$. Then the regular cover $\Gamma'$ of $\Gamma$ has vertices $G \times V$ and edges the collection of all pairs $(g,v)(g(l(v,w)),w)$ for each oriented edge $(v,w)$ of $\Gamma$ and each $g \in G$. 

We now show that $\Gamma'$ is also the graph of fibers of a ts of degree 2. Let $T$ be the 0-simple semigroup $\mathcal{M}^0(G,V \cup E,V,C)$, where $C$ is the $V \times (V \cup E)$ matrix over $G$
defined as follows. The $V \times V$ submatrix of $C$ is the identity matrix on $V$. If $(v,w)$ is an oriented edge of $\Gamma$, then in column $vw$ of $C$, 
$C(v,vw)=1$, $C(w,vw)=l(v,w)^{-1}$ and all other entries are 0. It is not difficult to check that $C$ has no distinct rows (columns) that are left (right) multiples of each other by an element of $G$. Thus $S$ is a $GM$ semigroup and by semi-local theory \cite[Section 4.6]{qtheory} we obtain a ts $(G \times V,T)$. 

\bl

With definitions as above, $(G \times V,T)$ is a ts of degree at most 2 and its fiber graph is the regular cover $\Gamma'$ of $\Gamma$.

\el 

\proof 
Since we assume that the degree of $(V,S)$ is 2, $E$ is non-empty. Let $(a,h,b) \in T$, where $a \in V \cup E, h \in G, b \in V$. Let $(g,u) \in G \times V$. If $a \in V$, then $(g,u)(a,h,b)$ is defined if and only if $u = a$. In this case, $(g,u)(a,h,b)=(gh,b)$ and it follows that the degree of $(a,h,b)$ is 1. Let $a = vw \in E$, with $(v,w) \in E^{+}$. Then $(g,u)(a,h,b)$ is defined if and only if $u=v$ or $u=w$ and we have $(g,v)(a,h,b) = (gh,b)$ and $(g,w)(a,h,b)=(g(l(v,w))^{-1}h,b)$. It follows that the fibers of $(a,h,b)$ are all sets of the form $\{(g,v),(g(l(v,w)),w)\}$. This both proves that the degree of  $(G \times V,T)$ is 2 and that its graph of fibers is 
$\Gamma'$. \qed

We now look at more detail at the graph of fibers of a $GM$ semigroup of degree at most 2 and the relationship to that of its $RLM$ image. Let $S$ be a $GM$ semigroup and let $I(S) \approx \mathcal{M}^{0}(G,A,B,C)$ be its distinguished 0-minimal ideal. Then we obtain a ts $(G \times B,S)$ which for simplicity's sake, we identify with $S$ throughout this discussion. Semi-local theory \cite[Section 4.6]{qtheory} tells us that $S$ induces a transformation semigroup $(B,RLM(S))$ and 
that $S$ embeds into the wreath product $G \wr (B,RLM(S))$. It follows that the complexity of $S$ is less than or equal to 1 plus the complexity of $RLM(S)$. Deciding whether the complexity of $S$ is equal or greater by 1 than that of $RLM(S)$ is the central problem of Krohn-Rhodes complexity theory as the computation of the complexity of an arbitrary semigroup can be reduced to this question \cite[Chapter 4]{qtheory}.

The following Lemma is useful for us.

\bl \label{wreathdeg}

Let $X=(Q,S)$ be a ts and $G$ a group. Then deg($G \wr X$)=deg($X$).

\el

\proof
Let $(f,s)$ be a transformation of $G \wr X$ and let $Y = \{(g_{1},q_{1}),\ldots , (g_{k},q_{k})\} \subseteq G \times Q$ be such that $Y(f,s)=\{(g,q)\}$ for some 
$(g,q) \in (G \times Q)$. Then $q_{i}s =q, i=1, \ldots , k$ and $g_{i}(q_{i})f=g, i=1, \ldots , k$. It follows that if $q_{i}=q_{j}$, then $g_{i}=g_{j}$. Therefore, 
$|Y| \leq deg(X)$ and thus, deg($G \wr X) \leq$ deg($X$). Let $s \in S$ be such that deg($s$)=deg($X$) and let $Y \subseteq Q$ be a fiber of $s$ of cardinality
deg($s$). Let $f:Q \rightarrow G$ be the constant function to the identity element 1 of $G$. Then $(\{1\} \times Y)(f,s)$ is a fiber of $G \wr X$ and therefore
deg($X) \leq $ deg$(G \wr X)$ and the result follows. \qed

If $S$ is a $GM$ semigroup then the embedding of $S$ into $G \wr (B,RLM(S))$ immediately implies the following corollary.

\bc \label{GMdegree2}

Let $S$ be a $GM$ semigroup. Then deg$(G \times B, S)$ = deg$(B,RLM(S))$.

\ec

We will now give a necessary and sufficient condition that deg$(G \times B, S) \leq 2$ for a $GM$ semigroup $S$ with distinguished ideal $I(S) \approx \mathcal{M}^{0}(G,A,B,C)$. Note that $I(S)$ is also a $GM$ semigroup and thus we have the ts $(G \times B, I(S))$. The following easy Lemma is useful.

\bl \label{fiberGM}

Let $S$ be a $GM$ semigroup with distinguished ideal $I(S)$. If  $F \subseteq (G \times B)$ is a fiber of some element of $S$ then it is a subset of a fiber of some element of $I(S)$. Therefore deg$(G \times B, S)$ = deg$(G \times B, I(S))$

\el

\proof

Let $s \in S$ and $(g,b) \in (G \times B)$ be such that $F=(g,b)s^{-1}$.  $(g,b)$ belongs to the distinguished $\mathcal{R}$-class of $I(S)$ and therefore, there is an idempotent $e$ in the $\mathcal{L}$-class of $(g,b)$ and thus $(g,b)e=(g,b)$. Therefore, $F=(g,b)s^{-1} \subseteq (g,b)e^{-1}s^{-1} =(g,b)(se)^{-1}$ and thus $F$ is contained in a fiber of $se \in I(S)$. Since $I(S) \subseteq S$ it easily follows that deg$(G \times B, S)$ = deg$(G \times B, I(S))$. \qed

\bc \label{deg2fiber}

Let $(G \times B,S)$ be a $GM$ semigroup of degree 2 and let $F$ be a fiber of size 2 of $X$. Then there is an idempotent $e \in I(S)$ such that $F$ is a fiber of $e$.
The same assertion holds for $(B,RLM(S))$.

\ec
\proof

By Lemma \ref{fiberGM}, $F$ is contained in a fiber of some element $s \in I(S)$. Since maximal fibers of $X$ have 2 elements, it follows that $F$ is a fiber of $s$. Since $I(S)$ is a regular semigroup, there is an idempotent $e \in I(S)$ and an element $x \in I(S)$ such that $sxs=es=s$. It easily follows that the partition 
$ker(s)=ker(e)$ and in particular, $F$ is a fiber of $e$. A similar proof proves the assertion for $(B,RLM(S))$. \qed

The next corollary follows immediately.

\bc \label{fiberideal}

Let $S$ be a $GM$ semigroup with distinguished ideal $I(S)$. Then the fiber graph of $(G \times B,S)$ is equal to the fiber graph of $(G \times B, I(S))$. The same result holds for $RLM(S)$ relative to its distinguished ideal.

\ec

We can use the preceding corollary to completely describe the fibers of a $GM$ semigroup of degree 2 in terms of its $RLM$ image. This is related to the ``Tie-Your-Shoes'' Lemma \cite[4.14.29]{qtheory}.

\bp \label{GMfiberRLM}

Let $X=(G \times B,S)$ be a $GM$ semigroup of degree 2. Then $F = \{(g',b'),(g,b)\}$ is a fiber of size 2 if and only if $\{b',b\}$ is a fiber of size 2 of $RLM(S)$ and there is an $a \in A$, such that  $g'C(b',a)C(b,a)^{-1}=g$.

\ep

\proof

Assume that $F$ is a fiber of size 2 of $S$. By Corollary \ref{fiberGM}, we can assume that $F$ is the fiber of an idempotent $e \in I(S)$. As is well known, the image of an idempotent partial function is a system of distinct representatives of its kernel. Therefore we can assume that $(g,b) \in Im(e)$ and thus that $(g,b)e=(g,b)$. It follows that $e=(a,C(b,a)^{-1},b)$ for some $a \in A$. Since $(g',b')e=(g,b)$, it follows that $C(b',a)\neq 0$ and that $g'C(b',a)C(b,a)^{-1}=g$. 
Clearly, the image of $e$ in $RLM(S)$ shows that $\{b,b'\}$ is a fiber of size 2 of $RLM(S)$. This proves the assertion in one direction.

Assume that $\{b,b'\}$ is a fiber of size 2 of $RLM(S)$ and there is an $a \in A$, such that  $g'C(b',a)C(b,a)^{-1}=g$. By Corollary \ref{fiberGM}, there is an idempotent
$f \in RLM(I(S))$ such that (without loss of generality), $b'f=bf=b$. As is well known, there is an idempotent $e \in S$ whose image in $RLM(I(S))$ is equal to $f$. It follows that $e=(a,C(b,a)^{-1},b)$ and by the assumption, it follows that $(g',b')e=(g,b)e=(g,b)$. Therefore, $\{(g,b),(g',b')\}$ is a fiber of size 2. \qed

The next proposition relates the graph of fibers of a $GM$ semigroup $S$ of degree 2 to that of its $RLM$ image. It summarizes the discussion above. If $\Gamma$ is a multigraph then its {\em simplification} is the simple graph one obtains by removing all loops from $\Gamma$ and replacing all multiple edges by a single edge.

\bp \label{fibercover}

Let $S$ be a $GM$ semigroup of degree 2  and let $I(S) \approx \mathcal{M}^{0}(G,A,B,C)$ be  its distinguished 0-minimal ideal. Let $\Gamma=(B,E)$ be the graph of fibers of 
$(B,RLM(S))$ and let $\Gamma'$ be the graph of fibers of $(G \times B,S)$. Then the following holds.

\begin{itemize}

	\item[(i)] {$\Gamma'=(G \times B,E')$ where $E'=\{\{(g',b'),(g,b)\}|\{b',b\} \in E \text{ and } g'C(b',a)C(b,a)^{-1}=g\}$.}
	
	\item[(ii)]{$G$ acts freely on $\Gamma'$ by $h(g,b)=(hg,b), h\{(g',b'),(g,b)\}=\{(hg',b'),(hg,b)\}$.}
	
	\item[(iii)] {The gain graph corresponding to the regular covering $\Gamma' \rightarrow \Gamma'/G$ is the multigraph with vertices $B$ and labelled edges 
	$(b',C(b',a)C(b,a)^{-1},b)$, where $C(b,a),C(b',a)$ are not 0 and $(b'b)$ is an edge of $\Gamma$.}
	
	\item[(iv)]{The simplification of the regular covering $\Gamma'/G$ is isomorphic to $\Gamma$.}

\end{itemize}

\ep
	
\proof

Item (i) follows immediately from Proposition \ref{GMfiberRLM}. Items (ii) and (iii) are clear from the definitions and item (i). Item (iv) follows from items (i),(ii) and (iii) and the definitions. \qed
	
Another important gain graph associated to a 0-simple semigroup $S \approx \mathcal{M}^{0}(G,A,B,C)$ is its Graham-Houghton graph, $GH(S)$. This graph appeared independently in the papers  of Graham \cite{Graham} and Houghton
\cite{houghton} and plays an important role in understanding various subsemigroups of 0-simple semigroups and their topology. See \cite[Section 4.13]{qtheory} 
where it is called the incidence graph of a Rees matrix semigroup. $GH(S)$ has vertices $A \cup B$ (where we assume that $A$ and $B$ are disjoint) and edges 
$E=\{ba|C(b,a) \neq 0\}$. We take oriented edges $E^{+}=\{(b,a)|ba \in E\}$ and label $(b,a)$ with $C(b,a)$. We see then that the graph of fibers $\Gamma'$ of $S$ as in Proposition \ref{fibercover} is the derived graph of the gain graph that has edges the labeled paths of length 2 in $GH(S)$ that begin at some $b' \in B$.

Regular covers with group $G$ as defined in this section are the combinatorial analogue of principal $G$-bundles. In this analogy, if $G$ is a group and $\Gamma$ is a 
graph with vertices $V$ and edges $E$, then the trivial $G$-cover is the graph with vertices $G \times V$ and edges $G \times E$ together with the projections to $V$ and $E$ respectively. $\Gamma$ is recovered by the obvious free actions of $G$ on $G \times V$ and $G \times E$. The trivial cover corresponds to the gain graph 
over $\Gamma$ and $G$ all of whose labels are the identity element of $G$. 

Houghton \cite{houghton} defined a cohomology theory for gain graphs and applies it to the theory of 0-simple semigroups. See also \cite[Section 4.13]{qtheory}. What Houghton defines are known as switching functions in the theory of gain graphs. We give brief details. Let $\Gamma = (V,E)$ be a graph. We define the group $B(\Gamma,G)$ of all functions $\delta:V \rightarrow G$ with pointwise multiplication. Two gain graphs $\Gamma(l)$ and $\Gamma(l')$ over $G$ are cohomologous if there is a $\delta \in B(\Gamma,G)$ such that for all oriented edges $e = (v,w)$, 
$l'(e)= \delta(v)l(e)\delta(w)^{-1}$. Such $\delta$ are called switching functions in the theory of gain graphs. $\Gamma(l)$ is called $G$-acyclic, if the label on any cycle in $\Gamma$, that is the product in $G$ of the values on each edge of the cycle, is the identity element of $G$. Then $\Gamma(l)$ is $G$-acyclic if and only if 
$l$ is cohomologous to the trivial labeling of $G$ \cite{houghton} or Proposition 4.13.14 of \cite{qtheory}. It is not hard to see that if $\Gamma(l)$ is cohomologous to $\Gamma(l')$, then the corresponding regular covers of $\Gamma$ are isomorphic. Furthermore, cohomologous gain graphs give isomorphic Rees matrix semigroups. See 
Section 4.13 of \cite{qtheory}.

An important theorem combining work of \cite{Graham} and \cite{houghton}  states that if
$S \approx \mathcal{M}^{0}(G,A,B,C)$ is a 0-simple semigroup, then its Graham-Houghton graph is $G$-acyclic if and only if $S$ is isomorphic to a 
Rees matrix semigroup $\mathcal{M}^{0}(G,A,B,C')$ in which all the
elements of $C'$ are 0 or 1. This latter condition is equivalent to having the subsemigroup of $S$ generated by its idempotents having only trivial subgroups. See Proposition 4.13.30 of \cite{qtheory}.

We apply these ideas to the relationship between the graph of fibers of a $GM$ semigroup and its $RLM$ image. Regular covers with group $G$ as defined in this section are the combinatorial analogue of principal $G$-bundles \cite{MM-Sheaves}. In this analogy, if $G$ is a group and $\Gamma$ is a 
graph with vertices $V$ and edges $E$, then the trivial $G$-cover is the graph with vertices $G \times V$ and edges $G \times E$ together with the projections to $V$ and $E$ respectively. $\Gamma$ is recovered by the obvious free actions of $G$ on $G \times V$ and $G \times E$. The trivial cover corresponds to the  cohomology class of the gain graph
over $\Gamma$ and $G$ all of whose labels are the identity element of $G$. Combining this with the discussion above, we have the following lovely connection between the topology and algebra of $GM$ semigroups of degree 2 and their $RLM$ image.

\bt \label{trivcov}

Let $S$ be a $GM$ semigroup of degree 2. The following conditions are equivalent.

\begin{itemize}

	\item[(1)]{The graph of fibers $\Gamma$ of $S$ is the trivial cover of the graph of fibers of $RLM(S)$.}
	
	\item[(ii)] $I(S)$ is isomorphic to a Rees matrix semigroup whose structure matrix has entries in $\{0,1\}$.
	
	\item[(iii)] The idempotent generated subsemigroup of $I(S)$ has only trivial subgroups.
	
\end{itemize}

\et

\proof

By Corollary \ref{fiberideal} we can assume that $S=I(S)$ is a 0-simple semigroup. Assume (i). If $\Gamma$ is the trivial cover over the graph of fibers of $RLM(S)$, then by the aforementioned connection between $\Gamma$ and the Graham-Houghton graph $GH(S)$, it follows that $GH(S)$ is $G$-acyclic. Therefore by Graham's theorem, $S$ 
is isomorphic to a Rees matrix semigroup whose structure matrix has entries in $\{0,1\}$. Thus (ii) holds. 

Graham's theorem states that (ii) and (iii) are equivalent. So assume (ii). By Proposition \ref{fibercover} (i) it follows that the graph of fibers of $S$ is the trivial cover over $RLM(S)$. \qed

We give some examples of fiber graphs of $GM$ and $RLM$ semigroups.

\be

Let $S \approx \mathcal{M}^{0}(Z_{2},\{a_{1},a_{2}\},\{b_{1},b_{2}\},C)$ where $C = 
	\begin{bmatrix}
1 & 1 \\
1 & -1		
	\end{bmatrix}$. One sees that $RLM(S)$ is the two element right zero semigroup with elements the constant functions to $b_1$ and $b_2$. The Graham-Houghton graph 
$GH(S)$ has underlying graph the complete bipartite graph $K(2,2)$ all labels are 1 except for the label from $b_2$ to $a_2$ which is -1. The graph of fibers of $S$ is the 4-cycle: $(1,b_{1})$---$(1,b_{2})$---$(-1,b_{1})$---$(-1,b_{2})$---$(1,b_{1})$. The quotient of this graph by the free action of $Z_2$ is the multigraph with two 
vertices and two edges between them. The fiber graph of $RLM(S)$ is the simplification of this graph, namely the path of length 1 with vertices $b_1$ and $b_2$.

\ee

We note that if we adjoin the identity element and the permutation that sends $(z,b_{i})$ to $(z,b_{i+1 (mod 2)}), z = \pm 1, i=1,2$ we obtain the semigroup discussed 
in Example \ref{small2}. The fact that the graph of fibers of this semigroup is a 4-cycle explains why this semigroup embeds into the monoid of continuous maps on a 4-cycle.
 
\be 

Let $S \approx \mathcal{M}^{0}(Z_{2},\{a_{1},a_{2}, a_{3},a_{4}\},\{b_{1},b_{2},b_{3},b_{4}\},C)$ where $C = 
	\begin{bmatrix}
1 & 1 & 0 & 0\\
0 & 1 & 1 & 0\\
0 & 0 & 1	& 1\\
1 & 0 & 0 & 1	
	\end{bmatrix}$. Then $RLM(S) \approx \mathcal{M}^{0}(\{a_{1},a_{2}, a_{3},a_{4}\},1,\{b_{1},b_{2},b_{3},b_{4}\},C)$ and its graph of fibers is the 4 cycle 
	$b_{i}$---$b_{i+1}, i(mod 4)$. By Theorem \ref{trivcov}, the graph of fibers of $S$ is the trivial cover of the 4-cycle relative to $Z_2$, that is, two disjoint copies of a 4-cycle. $GH(S)$ is the 8-cycle $b_{1}a_{1}b_{2}a_{2}b_{3}a_{3}b_{4}a_{4}b_{1}$ with all edges labeled by 1.
	
Now let $T \approx \mathcal{M}^{0}(Z_{2},\{a_{1},a_{2}, a_{3},a_{4}\},\{b_{1},b_{2},b_{3},b_{4}\},C)$ where $C = 
	\begin{bmatrix}
1 & 1 & 0 & 0\\
0 & 1 & 1 & 0\\
0 & 0 & 1	& 1\\
1 & 0 & 0 & -1	
	\end{bmatrix}$. Then $RLM(T)\approx RLM(S)$. The graph of fibers of $T$ is the 8-cycle 
	$(1,b_{1})$---$(1,b_{2})$---$(1,b_{3})$---$(1,b_{4})$---$(-1,b_{1})$---$(-1,b_{2})$---$(-1,b_{3})$---$(-1,b_{4})$---$(1,b_{1})$. $GH(S)$ is the gain graph,
	$b_{1}\stackrel{1}{\rightarrow}a_{2}\stackrel{1}{\rightarrow}b_{2}\stackrel{1}{\rightarrow}a_{3}\stackrel{1}{\rightarrow}b_{3}\stackrel{1}{\rightarrow}a_{4}\stackrel{-1}{\rightarrow}b_{4}\stackrel{1}{\rightarrow}a_{1}\stackrel{1}{\rightarrow}b_{1}$.

\ee
	
The preceding two examples show how both regular covers of degree 2 over a 4-cycle arise semigroup theoretically.



The results of this paper can easily be generalized to ts of arbitrary degree $k>2$. The fibers then have the structure of a hypergraph or an incidence system. Most of the basic theorems in this paper on ts of degree 2 can be suitably generalized to the case of arbitrary degree. We note that the complexity of a ts of degree $k$ is at most $k$ \cite{Tilsonnumber}.

As a special case, we have the problems considered in this paper for degree 2 $GM$ semigroups. As mentioned previously, every such semigroup has complexity 
either 1 or 2. This leads to the following problem.

\begin{Problem}

Let $S$ be a $GM$ semigroup of degree 2. Is there an algorithm to decide if the complexity of $S$ is 1?

\end{Problem}

In the continuation of this paper \cite{deg2part2}, we will give a number of illuminating examples of $GM$ semigroups of degree 2 and study in depth the problem above. In particular, we study the problem of whether the lower bound described in \cite{Trans} is an upper bound as well in the case of ts of degree 2.

\bibliography{stubib}

\def\malce{\mathbin{\hbox{$\bigcirc$\rlap{\kern-7.75pt\raise0,50pt\hbox{${\tt
  m}$}}}}}\def\cprime{$'$} \def\cprime{$'$} \def\cprime{$'$} \def\cprime{$'$}
  \def\cprime{$'$} \def\cprime{$'$} \def\cprime{$'$} \def\cprime{$'$}
\begin{thebibliography}{10}

\bibitem{projcont}
J.~H. Dinitz and S.~W. Margolis.
\newblock Continuous maps in finite projective space.
\newblock In {\em Proceedings of the thirteenth Southeastern conference on
  combinatorics, graph theory and computing (Boca Raton, Fla., 1982)},
  volume~35, pages 239--244, 1982.

\bibitem{bibdcont}
J.~H. Dinitz and S.~W. Margolis.
\newblock Continuous maps on block designs.
\newblock {\em Ars Combin.}, 14:21--45, 1982.

\bibitem{Eilenberg}
S.~Eilenberg.
\newblock {\em Automata, languages, and machines. {V}ol. {B}}.
\newblock Academic Press, New York, 1976.
\newblock With two chapters (``Depth decomposition theorem'' and ``Complexity
  of semigroups and morphisms'') by Bret Tilson, Pure and Applied Mathematics,
  Vol. 59.

\bibitem{Graham}
R.~L. Graham.
\newblock On finite {$0$}-simple semigroups and graph theory.
\newblock {\em Math. Systems Theory}, 2:325--339, 1968.

\bibitem{topgraph}
J.~L. Gross and T.~W. Tucker.
\newblock {\em Topological Graph Theory}.
\newblock Dover, New York, 2001.

\bibitem{graphmorphs}
P.~Hell and J.~N$\check{\text{e}}$set$\check{\text{r}}$il.
\newblock {\em Graphs and Homomorphisms}.
\newblock Oxford Lecture Series in Mathematics and its Applications. Oxford
  University Press, 2004.

\bibitem{Trans}
K.~Henckell, J.~Rhodes, and B.~Steinberg.
\newblock An effective lower bound for the complexity of finite semigroups and
  automata.
\newblock {\em Trans. AMS}, 364(4):1815--1857, 2012.

\bibitem{houghton}
C.~H. Houghton.
\newblock Completely {$0$}-simple semigroups and their associated graphs and
  groups.
\newblock {\em Semigroup Forum}, 14(1):41--67, 1977.

\bibitem{Arbib}
K.~Krohn, J.~Rhodes, and B.~Tilson.
\newblock {\em Algebraic Theory of Machines, Languages, and Semigroups}.
\newblock Edited by Michael A. Arbib. With a major contribution by Kenneth
  Krohn and John L. Rhodes. Academic Press, New York, 1968.
\newblock Chapters 1, 5--9.

\bibitem{Lawson}
M.~V. Lawson.
\newblock {\em Inverse semigroups}.
\newblock World Scientific Publishing Co. Inc., River Edge, NJ, 1998.
\newblock The theory of partial symmetries.

\bibitem{LyndonandSchupp}
R.~C. Lyndon and P.~E. Schupp.
\newblock {\em Combinatorial group theory}.
\newblock Classics in Mathematics. Springer-Verlag, Berlin, 2001.
\newblock Reprint of the 1977 edition.

\bibitem{MM-Sheaves}
S.~{Mac Lane} and I.~Moerdijk.
\newblock {\em Sheaves in geometry and logic}.
\newblock Universitext. Springer-Verlag, New York, 1994.
\newblock A first introduction to topos theory, Corrected reprint of the 1992
  edition.

\bibitem{deg2part2}
S.~Margolis and J.~Rhodes.
\newblock Degree 2 transformation semigroups: {C}omplexity and examples.
\newblock In preparation.

\bibitem{AmigoWilson}
S.~Margolis, J.~Rhodes, and P.~Silva.
\newblock On the {W}ilson monoid of a pairwise balanced design.
\newblock {\em J. Alg. Comb.}, To appear, 2021.

\bibitem{Tilsonnumber}
S.~W. Margolis.
\newblock {$k$}-transformation semigroups and a conjecture of {T}ilson.
\newblock {\em J. Pure Appl. Algebra}, 17(3):313--322, 1980.

\bibitem{bibdtrans}
S.~W. Margolis and J.~H. Dinitz.
\newblock Translational hulls and block designs.
\newblock {\em Semigroup Forum}, 27(1-4):247--263, 1983.

\bibitem{TilsonMargolis}
S.~W. Margolis and B.~Tilson.
\newblock An upper bound for the complexity of transformation semigroups.
\newblock {\em J. Algebra}, 73(2):518--537, 1981.

\bibitem{PetrichInv}
M.~Petrich.
\newblock {\em Inverse semigroups}.
\newblock Pure and Applied Mathematics (New York). John Wiley \& Sons Inc., New
  York, 1984.
\newblock A Wiley-Interscience Publication.

\bibitem{qtheory}
J.~Rhodes and B.~Steinberg.
\newblock {\em The {$q$}-theory of Finite Semigroups}.
\newblock Springer Monographs in Mathematics. Springer, New York, 2009.

\bibitem{Serre}
J.-P. Serre.
\newblock {\em Trees}.
\newblock Springer Monographs in Mathematics. Springer-Verlag, Berlin, 2003.
\newblock Translated from the French original by John Stillwell, Corrected 2nd
  printing of the 1980 English translation.

\bibitem{Steinpartial}
I.~Stein.
\newblock The representation theory of the monoid of all partial functions on a
  set and related monoids as {EI}-category algebras.
\newblock {\em Journal of Algebra}, 450, no. 15:549--569, 2016.

\end{thebibliography}
\bibliographystyle{abbrv}

\vspace{1cm}

{\sc Stuart Margolis, Department of Mathematics, Bar Ilan University,
  52900 Ramat Gan, Israel}

{\em E-mail address:} margolis@math.biu.ac.il

\bigskip

{\sc John Rhodes, Department of Mathematics, University of California,
  Berkeley, California 94720, U.S.A.}

{\em E-mail addresses}: rhodes@math.berkeley.edu, BlvdBastille@gmail.com

\bigskip



\end{document}